\input amstex
\documentstyle{amsppt}
\magnification=1150
\nopagenumbers
\centerline{\bf A Grazing Gaussian Beam}
\vskip.2in\centerline{James Ralston and Neelesh Tiruviluama}
\vskip.2in
\noindent {\bf Introduction}
\vskip.1in
In high frequency solutions of wave equations there is sharp propagation along ray paths, the projections into physical space of the bicharacteristic curves associated with the wave equation. Gaussian beams are asymptotic solutions concentrated on a single ray path in the high frequency limit. The interaction of Gaussian beams with boundaries has been studied in several settings. When a ray path intersects a boundary non-tangentially, the reflection of the beam is completely explicit, and it follows the reflected ray path. In the case of a boundary which is convex with respect to ray paths a suitably modified beam will follow the induced ray path on the boundary ([R]). Also, in the case of a cusped ray path, where the beam {\it Ansatz} fails, an asymptotic solution which is initially a beam will reappear as a beam after passing through the cusp ([Ti]).

A case which is not covered in these examples arises when the ray path is tangent to a concave boundary. This is usually called a \lq\lq grazing" ray. The propagation of singularities in this situation has been well understood since the 1980's: ordinary singularities simply continue to follow the ray path, [T]\footnote {Singularities in analyticity propagate in a more complicated fashion, [Sj] and [L]}. The question that we investigate here is how a beam following a grazing ray will interact with the boundary. Given the behavior of  beams in the cases above, it is likely that the asymptotic solution will again take the form of a gaussian beam outside a neighborhood in space-time of the point where the ray grazes the boundary. We have not been able to show that. However, in the example discussed here we are able to compute the amplitude of the asymptotic solution on the central ray and see the influence of grazing the boundary.

Since we are only studying an example, there are no theorems in this paper. Our main result is formula (23), combined with formula (24). These show that the amplitude of the solution on the central ray after grazing the boundary is very close to one half of what it was before hitting the boundary.
\vskip.2in
\noindent {\bf \S1.  Friedlander's Example}
\vskip.1in
F.G. Friedlander [F] observed that the  boundary value problem
$$Pu=(1+x)u_{tt}-u_{xx}-u_{yy}=0\hbox{ in }\{x>0\}\times \Bbb R_y\times \Bbb R_t,\ u(0,y,t)=f(y,t)\in C_c^\infty(\Bbb R^2),$$
where $u=0$ for $t<<0$, has the explicit solution
$$u(x,y,t)=(2\pi)^{-2}\int_{\Bbb R^2}e^{i(y\eta+t\tau)}{Ai(\zeta(x,\eta,\tau))\over Ai(\zeta(0,\eta,\tau))}\hat f(\eta,\tau)d\eta d\tau,\eqno{(1)}$$
where $Ai(z)$ is the standard Airy function, defined by the improper Riemann integral
$$Ai(x)={1\over 2\pi}\int_{\Bbb R}e^{i(x\xi+\xi^3/3)}d\xi$$
for $x\in \Bbb R$ and  extended to an entire function. The function $u(x,y,t)$ in (1) solves Friedlander's  wave equation when
$\zeta(x,\eta,\tau)= \beta(1+x-\eta^2/\tau^2)$, where $\beta^3=-\tau^2$. The choice of $\beta$ is determined by the requirement that $u=0$ for $t<<0$ and hence $\hat u$ is bounded analytic in Im$\{\tau\}<0$. The asymptotic formula in the sector $|\hbox{arg}(z)|<\pi-\delta$, $\delta>0$,
$$Ai(z)={1\over 2\sqrt \pi z^{1/4}}(1+O(z^{-3/2}))\exp(-2z^{3/2}/3)\eqno{(2)}$$
as $|z|\to\infty$, shows that when $x>0$ the quotient $Ai(\zeta(x,\eta,\tau))/ Ai(\zeta(0,\eta,\tau))$ will extend to a bounded function on Im$\{\tau\}<0$ if Re$\{\beta^{3/2}\}\geq 0$.
This is satisfied by taking $\beta=\tau^{2/3}\exp(\pi i/3)$ when $\tau>0$ and $\beta=\tau^{2/3}\exp(-\pi i/3)$ when $\tau<0$.

Friedlander's example is important in the study of propagation of singularities and high frequency asymptotics for wave equations because it contains the fundamental solution for a boundary value problem with \lq\lq grazing (aka glancing) rays". In the next section we will construct a gaussian beam that \lq\lq grazes" the boundary $x=0$.
\vskip.1in
\noindent {\bf \S2. A Gaussian Beam}
\vskip.1in
For simplicity we replace $P$ by $P/2$. Hence the symbol of $P$ becomes  $(\xi^2+\eta^2-(1+x)\tau^2)/2$, and its bicharateristic curves are solutions of
$$\dot x=\xi\qquad\dot y=\eta\qquad \dot t=-(1+x)\tau\qquad \dot \xi=\tau^2/2\qquad \dot \eta=0\qquad \dot \tau=0.$$
The gaussian beam that we will use here is constructed from the  bicharacteristic curve $(x(y),y,t(y),\xi(y),\eta(y),\tau(y))$, where
$$x(y)=y^2/4\qquad t(y)=y+y^3/12\qquad \xi(y)=y/2\qquad \eta(y)=1\qquad \tau(y)=-1.$$
Note that $(1+x(y))\tau^2(y)-\xi^2(y)-\eta^2(y)\equiv 0$, so this is a null bicharacteristic.
\vskip.1in
The gaussian beam we use here has the standard (first order) form
$v(x,y,t;k)=a(x,y,t)e^{ik\psi(x,y,t)}$ with $k>0$. We use the form $\psi_y=((1+x)(\psi_t)^2-(\psi_x)^2)^{1/2}$ for the eichonal equation, and require that the eichonal vanish to second order on the ray path $\gamma$ traced by $(x(y),y,t(y))$. This leads to the phase
$$\psi(x,y,t)=(x-x(y))\xi(y)+(t-t(y))\tau(y)+{1\over 2}(x-x(y),t-t(y))\cdot M(y)(x-x(y),t-t(y))$$
where $M(y)=W(y)(V(y))^{-1}$ with
$$W(y)=\left(\matrix i&-iy\\0&i\endmatrix\right)\hbox{ and } V(y)=\left(\matrix 1+iy&-i{y^2\over 2}\\y+i{y^2\over 2}&1-iy-i{y^3\over 4}\endmatrix\right).$$
That leads to
$$M(y)={i\over D}\left(\matrix 1-iy+y^2+i{y^3\over 4}&-y-i{y^2\over 2}\\ -y-i{y^2\over 2}&1+iy\endmatrix\right), \hbox{ where }D=1+y^2+i{y^3\over 4}.\eqno{(3)}$$
We choose the amplitude $a(y)=(\det V(y))^{-1/2}=D^{-1/2}$ which makes  $a(0)=1$. (See Appendix I for details of this computation.)
\vskip.2in
\noindent {\bf \S3. Asymptotic Solution of the Boundary Value Problem}
\vskip.1in
We want to solve the Dirichlet problem for $(1+x)u_{tt}=u_{xx}+u_{yy}$ on $x=0$ with $u(x,y,t)$ close to the gaussian beam $v(x,y,t,k)$ constructed in \S2 when $t<<0$. To do this precisely we would need to begin with $u(x,y,t,k)$ equal to the exact solution in $\Bbb R_{x,y}^2\times\Bbb R_t$ with $(u(x,y,0,k),u_t(x,y,0,k))=(v(x,y,0,k),v_t(x,y,0,k))$, and then use Friedlander's fundamental solution to find an exact solution $w(x,y,t,k)$ in $x>0$  which agrees with $u(x,y,t,k)$ on $x=0$. Then the solution that we want would be $u(x,y,t,k)-w(x,y,t,k)$. Since we are only interested in the behavior of the solutions when $k$ is large, we will replace the exact solution $u(x,y,t,k)$ by the gaussian beam $v(x,y,t,k)$. Hence
$$\hat f(\eta,\tau)=\int_{\Bbb R^2}e^{-i(y\eta+t\tau)}a(y)e^{ik\psi(0,y,t)}dydt.$$
 Note that for $k$ large, $v(0,y,t,k)$ is negligibly small outside $\{|y|\leq k^{-1/4+\epsilon}\}$ for any $\epsilon>0$. In evaluating Friedlander's fundamental solution on $v(0,y,t,k)$ we will replace $M(y)$ by $M(0)=iI$ and $a(y)$ by $a(0)=1$.
This amounts to computing the correction from reflection on $x=0$ for a \lq\lq beam-like" incoming wave in $t<<0$.
However, if one uses $M(y)$ and $a(y)$ from (3) in evaluating $v(0,y,t,k)$, the leading order reflected wave on the central ray, $w(x,2x^{1/2},2x^{1/2}+{2\over 3}x^{3/2},k)$ in (23), will be unchanged (see Appendix III).
\vskip.1in

Since the amplitude $a(0,y,t)$ does not depend on $t$, the integration in $t$ can be done explicitly. We substitute $s=t-t(y)$ for $t$, and change $y$ to $z$ in anticipation of substituting into (1). This yields
$$\hat f(\eta,\tau)=({2\pi\over k})^{1/2}\int_{\Bbb R}e^{-iz\eta}\exp[-i\tau (z+{z^3\over 12})-ik{z^3\over 8} -k{z^4\over 32}-{(\tau+k)^2\over 2k}]dz.$$
Substituting that into (1) and stretching variables $\eta=k\mu$, $\tau =k\nu$, we have
$$w(x,y,t)=({k\over 2\pi})^{3/2}\int_{\Bbb R^3}e^{ik\Phi(t,y,z,\mu,\nu)}{Ai(\zeta(x,k\mu,k\nu))\over Ai(\zeta(0,k\mu,k\nu))}dz d\mu d\nu, \hbox{ where }$$
$$\Phi=(y-z)\mu+t\nu  -\nu(z+{z^3\over 12})-{z^3\over 8}+i[{z^4\over 32}+{(\nu +1)^2\over 2}].\eqno{(4)}$$

We will be computing asymptotics in the region $x>0$.  Hence, $\zeta(x,k\mu,k\nu)\to \infty$ as $k\to\infty$, and we can use the asymptotic form (2). That gives
$$Ai(\zeta(x,k\mu,k\nu))={1\over 2\sqrt \pi}(\zeta(x,k\mu,k\nu))^{-1/4}(1+O(k^{-1}))\exp(-{2\over 3}\zeta(x,k\mu,k\nu)^{3/2}).\eqno{(5)}$$
Since the only significant contributions will come from $\nu$ close to $-1$, we use
 $$\zeta(x,k\mu,k\nu) =(\nu k)^{2/3}e^{-i\pi/3}(1+x-{\mu^2\over \nu^2}).$$
Also in $\zeta(x,k\mu,k\nu)^{3/2}$ we take $(\nu^{2/3})^{3/2}=|\nu|$, since the asymptotic formulas come from analytic continuation from the positive real axis.

For the factor $(Ai(\zeta(0,k\mu,k\nu))^{-1}$ the asymptotic form in (2) is not helpful, since $\zeta(0,k\mu,k\nu)=0$ when $\mu^2=\nu^2$. Let $\omega=\exp(2\pi i/3)$, and set $A(z)=Ai(\omega z)$. Following an observation in Lemma 9.9 in Melrose-Taylor [MT], we use the constant Wronskian, $W(z)=A(z)Ai^\prime(z)-A^\prime(z)Ai(z)$, to write
$${W(0)\over A(z)}=Ai^\prime(z)-{A^\prime(z)\over A(z)}Ai(z)={1\over 2\pi}\int_{\Bbb R}(i\xi-{A^\prime(z)\over A(z)})e^{i(z\xi+\xi^3/3)}d\xi.$$
Note, crucially, that when $\omega z=\zeta(0,k\mu,k\nu)$, $z=-(\nu k)^{2/3}(1-{\mu^2\over \nu^2})$. Letting $\xi=k^{1/3}T$
$$(Ai(\zeta(0,k\mu,k\nu))^{-1}=\hskip2in$$$${1\over 2\pi W(0)}\int_{\Bbb R}(ik^{1/3}T-e^{2\pi i/3}{Ai^\prime(\zeta(0,k\mu,k\nu))\over Ai(\zeta(0,k\mu,k\nu))})e^{ik[-\nu^{2/3}(1-{\mu^2\over \nu^2})T+T^3/3]}k^{1/3}  dT \eqno{(6)}$$
This gives the final asymptotic form for the reflected solution
 $$w(x,y,t.k)=\int_{\Bbb R^4}Z(k,x,\mu,\nu,T)e^{ik\Phi(t,x,y,z,\mu,\nu,T)}dTdz d\mu d\nu,\eqno{(7)}$$
where the phase $\Phi$ is given  by
$$\Phi(t,x,y,z,\mu,\nu,T)=(y-z)\mu+t\nu  -\nu(z+{z^3\over 12})-{z^3\over 8}+i[{z^4\over 32}+{(\nu +1)^2\over 2}]$$
$$+{2|\nu|\over 3}(1+x-{\mu^2\over \nu^2})^{3/2} -\nu^{2/3}(1-{\mu^2\over \nu^2})T+{T^3\over 3}\eqno{(8)}$$
and the amplitude $Z$ is given by
$$Z(k,x,\mu,\nu,T)\hskip3in$$$$={k^{11/6}\over \sqrt{ 2}(2\pi)^3W(0)}(ik^{1/3}T-e^{2\pi i/3}{Ai^\prime(\zeta(0,k\mu,k\nu))\over Ai(\zeta(0,k\mu,k\nu))})(\zeta(x,k\mu,k\nu))^{-1/4}(1+O(k^{-1})).\eqno{(9)}$$
\vskip.2in
\noindent {\bf \S 4. Stationary Phase Reduction}
\vskip.2in
In this section we remove two of the integrations in (7) by applying stationary phase. However, it is convenient to replace the variables $(\mu,\nu,T)$ by $(s,\nu, T)$ where $s=\mu+\nu$, and use stationary phase in $(s,T)$. The partial derivatives of the phase with respect to these variables are
$$\Phi_s=y-z-\nu(x+s({2\nu-s\over \nu^2}))^{1/2}({2\nu-2s\over \nu^2})-\nu^{-4/3}(2\nu-2s)T, \hbox{ and }$$
$$\Phi_T=T^2-\nu^{-4/3}s(2\nu-s),\eqno{(10)}$$
where we have assumed that $\nu<0$ in the computation of $\Phi_s$, since the only significant contributions to $w(x,y,t)$ will come from $\nu$ near -1.

The phase $\Phi$ will be stationary in $(s,T)$ when $\Phi_T=\Phi_s=0$. Note that we are only looking for stationary points where $x>0$.
One sees from (10) that $\Phi_T=\Phi_s=0$ on the surface $\Sigma =\{(x,y,z): 2\sqrt x=(y-z)\}$ when $(T,s)=(0,0)$.
At these points, the hessian determinant satisfies
$$ \hbox{ det}\left(\matrix \Phi_{ss}&\Phi_{sT}\\ \Phi_{Ts}&\Phi_{TT}\endmatrix\right)=-4\nu^{-2/3}.$$
Hence by the implicit function theorem, the phase is stationary on a set $(s,T)=(s(x,y,z,\nu),T(x,y,z,\nu))$, where the functions $s(x,y,z,\nu)$ and $T(x,y,z,\nu)$ are  defined on a neighborhood of $\Sigma\cap\{x>0\}$ and vanish on $\Sigma$. The structure of the stationary set is quite implicit in (10), and we found it easier to describe it using a change of variables related to the bicharacteristic flow for Friedlander's wave equation.
\vskip.2in
Using $y-z$ as the curve parameter in the solution of (A.1) with initial data $(x_0,t_0)=(0,0)$ one finds
$$x={\xi_0\over \eta}(y-z)+{\tau^2\over \eta^2}{(y-z)^2\over 4} \hbox{ and }t=-{\tau\over\eta}(y-z)-{\xi_0\tau\over \eta^2}{(y-z)^2\over 2}-{\tau^3\over \eta^3}{(y-z)^3\over 12} $$
with the constraint $\tau^2=\xi_0^2+\eta^2$ so that $(x,y,t)$ will lie on a null bicharacteristic.  A short computation shows  that,  using $\xi_0/\tau$ and $\eta/\tau$ as coordinates on the unit circle, $\Phi_T=0$ precisely when
$$(T,s)=(\nu^{1/3}{\xi_0\over \tau},\nu(1+{\eta\over \tau})).$$
The surface $\Sigma$ is composed of the projections of null bicharacteristics with data  $(\xi_0/\tau,\eta/\tau)=(0,-1)$. Hence we can cover the data of null bicharacteristics over a neighborhood of $\Sigma$ by restricting $\eta/\tau$ to be negative.

To study $\Phi^{sp}$, the phase evaluated on the stationary set, we introduce $r=\eta/\tau$. Then $\xi_0/\tau=\pm(1-r^2)^{1/2}$, and
 $$\Phi_s= y-z+2r[(x+(1-r^2))^{1/2}\pm(1-r^2)^{1/2}].\eqno{(11)} $$
$\Phi_s=0$ implies the quadratic equation for $r^2$
$$r^4(x^2+(y-z)^2)-r^2({x\over 2}+1)(y-z)^2+{(y-z)^4\over 16}=0.$$
We always assume $4+4x\geq(y-z)^2$ so that $r^2$ will be real.
 The relation of the roots of this equation to solutions of (11) is as follows. The root\footnote{The other root of the quadratic equation for $r^2$,$$r^2=\big[{x/2+1-\sqrt{1+x-(y-z)^2/4}\over 2(x^2+(y-z)^2)}\big](y-z)^2$$ also gives a solution to $0=\Phi_s= y-z+2r[(x+(1-r^2))^{1/2}+(1-r^2)^{1/2}]$ with $r<0$, but this does not contribute to the leading terms as $k\to\infty$. See \S 6.}
 $$r=-\big[{x/2+1+\sqrt{1+x-(y-z)^2/4}\over 2(x^2+(y-z)^2)}\big]^{1/2}(y-z)\eqno{(12)}$$
is a solution to
$$y-z+2r[(x+(1-r^2))^{1/2}-(1-r^2)^{1/2}]=0\hbox{ when }4x\geq (y-z)^2,$$
and a solution to
$$y-z+2r[(x+(1-r^2))^{1/2}+(1-r^2)^{1/2}]=0\hbox{ when }4x\leq (y-z)^2\leq 4x+4.$$
Thus at the stationary points (all square roots are positive)
$$x=\pm({1\over r^2}-1)^{1/2}(y-z)+{(y-z)^2\over 4r^2},\quad T=\mp\nu^{1/3}(1-r^2)^{1/2}\hbox{ and } s=\nu(1+r)  $$
with $r$ given by (12). In this formula and in what follows, the $+$ sign corresponds to $4x>(y-z)^2$ and the $-$ sign corresponds to $4x<(y-z)^2$. We will show (see (A.8)) that $r(x,y,y-2\sqrt x +p)=-1+p^2/8 +O(p^3)$. This implies that we can set $\pm (1-r^2)^{1/2}= p/2+O(p^2)$, and have the sign of $p$ correspond the correct choice of $\pm$ with respect to $(x,y)$. However, we will  continue to use $\pm (1-r^2)^{1/2}$ instead of $ p/2+O(p^2)$ until \S 6. In \S 6 we are using $y=2\sqrt x$ so that $p=z$
\vskip.1in
Using $r$ and $|\nu|=-\nu$, we have
$$\Phi^{sp}(t,x,y,\nu,z)\hskip3.5in \eqno{(13)}$$ $$=\nu[t-z-{z^3\over 12}+r(y-z)+{2\over 3}({y-z\over 2r}\mp(1-r^2)^{1/2})^3+{2\over 3}(\pm(1-r^2)^{1/2})^3]-{z^3\over 8}+{i\over 2}[(\nu+1)^2+{z^4\over 16}].\hskip2in$$
In deriving (13) one needs to be careful with signs, particularly in evaluating
$${2|\nu|\over 3}(1+x-{\mu^2\over \nu^2})^{3/2}=-{2\nu\over 3}(x+s{(2\nu-s)\over \nu^2})^{3/2}$$ $$=-{2\nu\over 3}
(\mp{1-r^2)^{1/2}\over r}(y-z)+{(y-z)^2\over 4r^2}+(1-r^2))^{3/2}$$
$$=-{2\nu\over 3}(\pm(1-r^2)^{1/2}-{y-z\over 2r})^3={2\nu\over 3}({y-z\over 2r}\mp(1-r^2)^{1/2})^3,$$
because $(1+x-{\mu^2\over \nu^2})^{3/2}$ was chosen to be positive, and $\pm(1-r^2)^{1/2}-{y-z\over 2r}$ is the positive square root of $\mp{1-r^2)^{1/2}\over r}(y-z)+{(y-z)^2\over 4r^2}+(1-r^2)$

One gets a further simplification by writing $\pm(1-r^2)^{1/2}$ as ${y-z\over 4r}-{rx\over y-z}$. Then combining the cubic terms in (13) gives the final form of $\Phi^{sp}$:
$$\Phi^{sp}=\nu[t-z-{z^3\over 12}+r(y-z)+{(y-z)^3\over48r^3}+{rx^2\over y-z}]-{z^3\over 8}+{i\over 2}[(\nu+1)^2+{z^4\over 16}]\eqno{(14)}$$
In these variables, on the stationary set
$$ J=_{def}\hbox{ det}\left(\matrix \Phi_{ss}&\Phi_{sT}\\ \Phi_{Ts}&\Phi_{TT}\endmatrix\right)=\hskip2in\eqno{(15)} $$$$4\nu^{-2/3}[(x+(1-r^2))^{-1/2}r^2(1-r^2)^{1/2}
-(x+(1-r^2))^{1/2}(1-r^2)^{1/2}+(1-r^2)-r^2]$$
Note $|J|=-J$, since $J$ is negative when $r=-1$.
\vskip.1in
The stationary phase expansion as given in Theorem 7.7.6 in [H] is adequate for our purposes here. The error in that expansion is bounded by derivatives of the amplitude in $(s,T)$, and we need more details on the amplitude $Z(k,x,s-\nu,\nu,T)$ in (9). Since the asymptotic expansion in (2) can be differentiated,
$${Ai^\prime\over Ai}(z)=-z^{1/2}-{1\over 4}z^{-1}+O(z^{-5/2})$$
in the sector $|\hbox{arg}(z)|<\pi-\delta$. Substituting $z=\zeta(0,k\mu,k\nu) =k^{2/3}\nu^{-4/3}e^{-i\pi/3}(2\nu s-s^2)$, one sees that the the error when one uses the first term in the stationary phase expansion, i.e.
$$Z(k,x,s(x,y,z,\nu),\nu, T(x,y,z,\nu))e^{ik\Phi^{sp}(x,y,z,t,\nu)},$$
 is order $k$. Since the stationary set is $(s,T)=(\nu(1+r),\mp\nu^{1/3}(1-r^2)^{1/2})$, the stationary phase expansion gives
$$w(x,y,t,k)=$$$$\int_{\Bbb R^2}({2\pi\over k})Z(k,x,\nu r, \nu,\mp\nu^{1/3}(1-r^2)^{1/2})+O(k))|J|^{-1/2}e^{ik\Phi_{sp}(t,x,y,z,\nu)}dzd\nu,\eqno{(16)}$$
where $Z$ is given by (9), $\Phi^{sp}$ is given by (14), and $J$ is from (15). Note that the signature of the hessian of $\Phi_{sp}$ is zero, because at $(s,T)=(0,0)$ this hessian has the form $\left(\matrix \cdot &\cdot\\
\cdot &0\endmatrix \right)$.
\vskip.2in
\noindent {\bf \S 5. Integration in $\nu$.}
\vskip.2in
At this point (7) has been reduced to
$$w(x,y,t,k)=\int_{\Bbb R^2}A((x,y,z,t,\nu,k)e^{ik(B+\nu C+i(\nu+1)^2/2)}d\nu dz, \eqno{(17)}$$
where
$$B(z)= -{z^3\over 8}+i{z^4\over 32}\hbox{ and }C(x,y,z,t)=t-z-{z^3\over 12}+r(y-z)+{(y-z)^3\over48r^3}+{rx^2\over y-z},$$
and $A(x,y,z,t,\nu,k)$ is the amplitude in (16).
The integral representation (17) is a gaussian integral in $\nu$. The contributions away from $\nu=-1$ decay exponentially as $k\to\infty$. The method of steepest descent applied to the integral in $\nu$ in (17) yields
$$\int_{\Bbb R}A(x,y,z,t,\nu,k)e^{ik(B+\nu C+i(\nu+1)^2/2)}d\nu=$$$$({2\pi\over k})^{1/2}A(x,y,z,t,-1+iC,k)e^{ik(B-C)-kC^2/2}(1+O(k^{-1})).\eqno{(18)}$$
\vskip.2in
\noindent \S 6. {\bf Discussion of the Final Integral.}
\vskip.2in
The ray $\gamma$, traced by $(x,y,t)=(y^2/4,y, y+y^3/12)$, is the projection of the null bicharacteristic for Friedlander's wave equation with initial data $(x,y,t,\xi,\eta, \tau)=(0,0,0,0,1,-1)$. Propagation of singularities  implies that $w(x,y,t,k)$ will be concentrated on this ray in the region $x>0,\ t>0$, i.e. when $y>0$.

 Since $w(x,y,t,k)$ is given by
$$({2\pi\over k})^{1/2}\int_{\Bbb R}A(x,y,t,k,-1+iC,z)e^{ik(B-C)-kC^2/2}(1+O(k^{-1}))dz,\eqno{(19)}$$
outside any neighborhood of $\{(x,y,z,t): C=0\}\cap\{z=0\}$ $w(x,y,t,k)$ will be $O(k^{-\infty})$.  Hence the asymptotic form of $w$ to leading order is determined by the Taylor series of $i(B-C)-C^2/2$ in $z$ about $z=0$ at points where $C(x,y,0,t)=0$.  The computations in Appendix II show that, on $\gamma$, $C(x,y,0,t)=\partial_zC(x,y,0,t)= \partial_z^2C(x,y,0,t)=0$.  Hence $C^2=O(z^6)$ on $\gamma$, and it will not contribute to the leading terms in the Taylor series of $i(B-C)-C^2/2$.
\vskip.1in
 We will conclude this investigation by computing the leading asymptotic form of $w$ on $\gamma$ as $k\to\infty$.
 From (12) one sees  that, when $x=(y-z)^2/4$, we have $r^2\equiv 1$ and, since $r$ is necessarily negative, this means $r=-1$. Further computation in Appendix II shows that $\partial_zr(x,2x^{1/2},0)=0$ and $\partial_z^2r(x,2x^{1/2},0)=1/4$. Combined with $\partial^3_zr(x,2x^{1/2},0)= (3x^{-1/2}-3x^{1/2})/8$ and
$\partial^4_zr(x,2x^{1/2},0)=15(x-1+x^{-1})/16$, this leads to
$\partial^2_z[B-C](x,2x^{1/2},0,t)=\partial^3_z[B-C](x,2x^{1/2},0,t) =0$ and
$$\partial^4_z[B-C](x,2x^{1/2},0,t)= -{3\over 8}(x^{1/2}-x^{-1/2})+3i/4 =_{def}24 a(x)$$ (see (A.12), (A.16) and (A.17)). Hence, using the parametrization of $\gamma$, $(x,y,t)=(x,2x^{1/2},2x^{1/2}(1+x/3))$,
$$[i(B-C)-C^2/2](x,2x^{1/2},z,t)=i(t-2x^{1/2}-{2\over 3}x^{3/2}+a(x)z^4(1+O(z))),\eqno{(20)}$$
where $1+O(z)$ denotes a power series with a constant term 1.
\vskip.1in
For the amplitude we recall that on the stationary set  when $y=2x^{1/2}$ $$(s,T)=(\nu(1+r),\mp\nu^{1/3}(1-r^2)^{1/2})=(-z ^2/8+O(z^3),z/2+O(z^2)).$$
Hence $\zeta(x,k\mu,k\nu)
=k^{2/3}e^{-i\pi/3}(x+(z^2/ 4)(1+O(z))=_{def}\xi(x,z,k).$
Substituting that into $Z$ in (9) we have
$$({2\pi\over k})^{1/2}A(x,2x^{1/2},t,k,-1+iC,z)=$$ $$
(2\pi)^{-3/2}{k^{1/6}e^{i\pi/12}\over 2^{1/2}W(0)}(ik^{1/3}{z\over 2}(1+O(z))-e^{2\pi i/3}{ Ai^\prime(\xi(0,z,k))\over Ai(\xi(0,z,k))})(x+O(z^2))^{-1/4}|J|^{-1/2}$$
and $J=-4(1+O(z))$ -- note that the $\zeta^{-1/4}$ in $Z$ introduces a factor of $k^{-1/6}$.  When we substitute these formulas into (19), we can make the change of variables $z=k^{-1/4}u$. This makes all of the $O(z)$ terms $O(k^{-1/4}u)$, and we are left with
$$w(x,2x^{1/2},2x^{1/2}+{2\over 3}x^{3/2},k)=$$$${ck^{1/6}\over x^{1/4}}\int_{\Bbb R}({i\over 2}k^{1/12}u -e^{2\pi i/3}{Ai^\prime\over Ai}(e^{-i\pi/3}{u^2\over 4}k^{1/6})e^{ia(x)u^4})k^{-1/4}du +O(k^{-1/2}) $$
$$={c\over x^{1/4}}\int_{\Bbb R} ({iu\over 2}-k^{-1/12}e^{2\pi i/3}{ Ai^\prime\over Ai}(e^{-i\pi/3}{u^2\over 4}k^{1/6}))e^{ia(x)u^4}du +O(k^{-1/2}),\eqno{(21)}$$
where
$$c=(4\pi)^{-3/2}{e^{i\pi/12}\over W(0)}.$$
Using the asymptotic behavior of $Ai^\prime(z)/Ai(z)$ and the dominated convergence theorem, we arrive at
$$w(x,2x^{1/2},2x^{1/2}+{2\over 3}x^{3/2},k)\sim {c\over 2x^{1/4}}\int_{\Bbb R}(iu+i|u|)e^{ia(x)u^4}du,\eqno{(22)}$$
as $k\to\infty$. The integral in (22) can be evaluated in the following way:
$$\int_{\Bbb R}ue^{-bu^4}du\qquad \hbox{ and } \qquad\int_0^\infty ue^{-bu^4}du$$
are analytic in $\hbox{Re}\{b\}>0$, and, setting $u=b^{-1/4}s$ for $b>0$, one sees that the first integral vanishes,  and
$$\int_0^\infty ue^{-bu^4}du=b^{-1/2}\int_0^\infty se^{-s^4}ds={1\over 4}b^{-1/2}\int_0^\infty t^{-1/2}e^{-t}dt={1\over 4}b^{-1/2}\Gamma({1\over 2}).$$
Both these relations hold in $\hbox{Re}\{b\}>0$ by analyticity. Using $b=-ia(x)$ and $W(0)=(e^{2\pi i/3}-1)/(2\pi\sqrt 3)$, formula (22) appears to be
$$w(x,2x^{1/2},2x^{1/2}+{2\over 3}x^{3/2},k)\sim {3\over 2(-3+3x-6ix^{1/2})^{1/2}(e^{2\pi i/3}-1)}e^{i\pi/3}$$$$={1\over 2}(1-x+2ix^{1/2})^{-1/2},\eqno{(23)}$$
as $k\to \infty$. The Gaussian beam from \S2 satisfies (see (3))
$$v(x,2x^{1/2},2x^{1/2}+{2\over 3}x^{3/2},k)=(1+4x+2ix^{3/2})^{-1/2},\eqno{(24
)}$$
and the  amplitude of the \lq\lq reflected" beam when  $x>0$ is
$$v(x,2x^{1/2},2x^{1/2}+{2\over 3}x^{3/2},k)-w(x,2x^{1/2},2x^{1/2}+{2\over 3}x^{3/2},k).$$
\vskip.1in
Continuing the computations in Appendix II, one can find the derivatives $C_x$ and $C_y$ on $4x=(y-z)^2$. This leads to
$$\partial_xw(x,y,t,k)\sim ikx^{1/2}w(x,y,t)\hbox{ and }\partial_yw(x,y,t,k)\sim ikw(x,y,t,k)\hbox{ and }$$
on $(x,y,t)=(x,2x^{1/2},2x^{1/2}+{2\over 3}x^{3/2})$. These results are consistent with $w$ being a gaussian beam concentrated on the ray path. One can also compute the higher derivatives $C_{xz} =0$, $C_{xzz}=3x^{-1/2}/4$, $C_{yz}=0$ and $C_{yzz}=1/4$ on $4x=(y-z)^2$. These suggest that the lower order terms in the expansion of $\partial_xw(x,y,t,k)$ and $\partial_yw(x,y,t,k)$ are {\it not} consistent with $w$ being a gaussian beam. However, we cannot reach that conclusion, since we only know the expansion of $w$ to leading order (see (21)).

\newpage

\centerline {\bf Appendix I}
 \vskip.2in
 The computation of $M(y)$ can be done fairly easily following the method in [R1]. For the symbol $\eta-((1+x)\tau^2-\xi^2)^{1/2}$ the bicharacteristic equations are
$$\dot x={\xi\over \eta}\qquad \dot y =1\qquad \dot t=-{\tau(1+x)\over\eta} \qquad \dot \xi={\tau^2\over 2\eta}\qquad \dot \eta=0\qquad \dot \tau=0. \eqno{(A.1)}$$
So we identify $y$ with the curve parameter, set $\eta=1$, and consider the reduced system for $(x,t,\xi,\tau)$ as functions of $y$.
We need the variation of solutions to (A.1) in this reduced form about the solution with the initial data $(x_0,t_0,\xi_0,\tau_0)=(0,0,0,-1)$.  The general solution to (A.1) is
$$x=x_0+\xi_0y+{\tau^2_0\over 4}y^2,\quad t=t_0-\tau_0((1+x_0)y+{\xi_0\over 2}y^2+{\tau^2_0\over 12}y^3)\quad \xi=\xi_0+{\tau^2_0\over 2}y.\quad \tau=\tau_0$$
So the variation is
$$\delta x=\delta x_0+(\delta \xi_0)y-{1\over 2}(\delta \tau_0)y^2\qquad \delta t=\delta t_0-(\delta\tau_0)(y+{y^3\over 12})+(\delta x_0)y+{\delta \xi_0\over 2}y^2-{\delta\tau_0\over 6}y^3$$
$$\delta \xi=\delta \xi_0-(\delta \tau_0)y\qquad \delta \tau=\delta \tau_0.$$
Letting $( v_{11},v_{12}, w_{11},w_{12})$ be the variation with $(\delta x_0,\delta t_0,\delta \xi_0,\delta \tau_0)=(1,0,i,0)$ and $( v_{21},v_{22}, w_{21},w_{22})$ be the variation with $(\delta x_0,\delta t_0,\delta \xi_0,\delta \tau_0)=(0,1,0,i)$, we have
$$M(y)=\left(\matrix w_{11}&w_{21}\\w_{21}&w_{22}\endmatrix\right)\left(\matrix v_{11}&v_{21}\\v_{12}&v_{22}\endmatrix\right)^{-1}=_{def}W(y)(V(y))^{-1}.$$
and we choose the amplitude $A(y)=(\det V(y))^{-1/2}$ with  $A(0)=1$.

For the amplitude recall that for a first order beam we need
$${d\over dy}(a(x(y),y,t(y))+{1\over 2}((1+x)\psi_{tt}-\psi_{xx}-\psi_{yy})|_\gamma a(x(y,y,t(y))=0$$
Using $\eta(x,\xi,\tau)=((1+x)\tau^2-\xi^2)^{1/2}$, and note that
$$\eta_\xi(x,\psi_x,\psi_t)=-{\psi_x\over \eta}\hbox{ and }\eta_\tau(x,\psi_x,\psi_t)={(1+x)\psi_t\over \eta}.$$
On $\gamma$ we have $\eta =\psi_y=1$ and $\psi_{yy}+\psi_{yx}\dot x+\psi_{yt}\dot t=0$. So
$$(1+x)\psi_{tt}-\psi_{xx}-\psi_{yy}=(\eta\eta_\tau)_t+(\eta\eta_\xi)_x-\psi_{yy}=(\eta_\xi)_x +(\eta _\tau)_t$$

Going back to the equations for $V$ and $W$ above one sees that on $\gamma$
$$(\eta_\xi)_x +(\eta_\tau)_t= \hbox{trace}({d V\over dy} V^{-1})={d\over dy}(\log\hbox{det} V).$$
Now the equation for the $a$ becomes
$${d\over dy}(\log a)=-{1\over 2}{d\over dy}\log(\hbox{det}V)$$
and $a$ is a constant multiple of $(\hbox{det}V)^{-1/2}$.
\vskip.2in
\centerline  {\bf Appendix II}
 \vskip.2in
 These are the computations that lead to the results in \S 6. Hence we need the  $z$-derivatives of $C$ (see (17)) when $z=0$.
\vskip.1in
The identity
$$0=r^4(x^2+(y-z)^2)-r^2({x\over 2}+1)(y-z)^2+{1\over 16}(y-z)^4, \eqno{(A.2)}$$
can  be differentiated with respect to $z$ to give
 $$0=4r^3(x^2+(y-z)^2)r_z-2r^4(y-z)-2r({x\over 2}+1)(y-z)^2r_z+2r^2({x\over 2}+1)(y-z)-{1\over 4}(y-z)^3.$$
 Evaluating that at $r=-1$ and $4x=(y-z)^2$, one sees that $r_z=0$ when $r=-1$ and $4x=(y-z)^2$.

 Taking the second derivative of (A.2) with respect to $z$
 $$0=(12r^2(r_z)^2+4r^3r_{zz})(x^2+(y-z)^2)+2r^4-16r^3(y-z)r_z$$$$ -2((r_z)^2+rr_{zz})({x\over 2}+1)(y-z)^2 +8r({x\over 2}+1)(y-z)r_z-2r^2({x\over 2}+1)+{3\over 4}(y-z)^2.\eqno{(A.3)}$$
 Evaluating that at $r=-1$, $r_z=0$ and $4x=(y-z)^2$ gives
 $$0=-4({1\over 16}(y-z)^4+(y-z)^2)r_{zz}+2+2({1\over 8}(y-z)^2+1)(y-z)^2r_{zz}-2({1\over 8}(y-z)^2+1)+{3\over 4}(y-z)^2\hskip1in$$
 which simplifies to
 $$0=-2(y-z)^2r_{zz}+{1\over 2}(y-z)^2.\eqno{(A.4)}$$
 Since $r=-1$ and $r_z=0$ are consequences of $4x=(y-z)^2$, (A.4) simply says
 $r_{zz}=1/4$ when $4x=(y-z)^2$.
 \vskip.1in
 Going back to (A.3) and differentiating with respect to $z$ once more
 $$0=(24r(r_z)^3+36r^2r_zr_{zz}+4r^3r_{zzz})(x^2+(y-z)^2)-2(12r^2(r_z)^2+4r^3r_{zz})(y-z)$$ $$-(48r^2(r_z)^2+16r^3r_{zz})(y-z)
 +24r^3r_z-2(3r_zr_{zz}+rr_{zzz})({x \over 2}+1)(y-z)^2$$$$+4((r_z)^2+rr_{zz})({x \over 2}+1)(y-z)+8((r_z)^2+rr_{zz})({x \over 2}+1)(y-z)-12rr_z({x \over 2}+1)-{3\over 2}(y-z).\eqno{(A.5)}$$

 Since $r=-1$, $r_z=0$ and  $r_{zz}=1/4$ when $4x=(y-z)^2$, and we are interested in the case $y-z\geq 0$, on $4x=(y-z)^2$ ((A.5) reduces to
 $$0=-4(x^2+4x)r_{zzz}+12\sqrt x+(4x^2+8x)r_{zzz}-6\sqrt x({x\over 2}+1)-3\sqrt x, $$
 which gives
 $$r_{zzz}={3\over 8}({1\over \sqrt x}-\sqrt x)\hbox{ when } y-z=2\sqrt x.\eqno{(A.6)}$$

We will need $r_{zzzz}$ but only when $y-z=2\sqrt x$. So differentiating (A.5) with respect to $z$ and evaluating on
$r=-1$, $r_z=0$ and  $r_{zz}=1/4$, we have
$$0=({9\over 4}-4r_{zzzz})(x^2+(y-z)^2)+32r_{zzz}(y-z)-12-16r_{zzz}({x\over 2}+1)(y-z)$$$$+(-{3\over 8}+2r_{zzzz})({x\over 2}+1)(y-z)^2+6({x\over 2}+1)+{3\over 2}. \eqno{(A.7)}$$
Substituting $r_{zzz}$ from (A.6) into (A.7) and solving for $r_{zzzz}$ we find
$$r_{zzzz}={15\over 16}(x-1+{1\over x})\hbox{ when }4x=(y-z)^2.$$
Hence collecting the results so far in this Appendix,
$$r(x,y,y-2\sqrt x+w)=-1+{w^2\over 8}+c_3(x){w^3\over 6}+c_4(x){w^4\over 24}+O(w^5),\eqno{(A.8)}$$
where $c_3(x)=(3/\sqrt x-3\sqrt x)/8$ and $c_4(x)=15(x-1+x^{-1})/16$.

 To arrive at (20) we need the first five terms in the Taylor series about $w=0$ for the phase $\Phi^{sp}(t,x,y,\nu,z)$ in (14) where $z=y-2\sqrt x+w$. Writing $\Phi^{sp}$ as
 $$\Phi^{sp}=\nu\phi(x,y,z)+t\nu-{\nu z^3\over 12}-{z^3\over 8}+{i\over 2}[(\nu+1)^2+{z^4\over 16}],$$
 we have from (14)
 $$\phi=-z+r(y-z)+{(y-z)^3\over 48r^3}+{rx^2\over y-z},\eqno{(A.9)}$$
 $$ \phi_z=-1+(y-z)r_z-r-{(y-z)^2\over 16r^3}+[-{(y-z)^3\over 16r^4}+{x^2\over y-z}]r_z+{rx^2\over (y-z)^2}, \eqno{(A.10)}$$
 $$\phi_{zz}=(y-z)r_{zz}-2r_z+2{(y-z)\over 16r^3}+6{(y-z)^2\over 16r^4}r_z-{(y-z)^3\over 16r^4}r_{zz}+{(y-z)^3\over 4r^5}(r_z)^2$$
 $$+{x^2\over (y-z)^2}r_z+{2rx^2\over (y-z)^3}+{x^2\over y-z}r_{zz}+{x^2\over (y-z)^2}r_z\eqno{(A.11})$$

When we evaluate (A.9), (A.10) and (A.11) at $4x=(y-z)^2
$ which implies $r=-1$, $ r_z=0$ and $r_{zz}=1/4$, we have
$$\phi=-y-{1\over 12}(y-z)^3, \qquad
\phi_z=0\quad\hbox{ and }\quad\phi_{zz}=0.\eqno{(A.12)}$$
These results may appear surprising but they just reflect the fact that $4x=(y-z)^2$ is the grazing set.

The next Taylor coefficient that we need is $\phi_{zzz}$. This, unfortunately, entails differentiating (A.11).
We have
$$\phi_{zzz}=-r_{zz}+(y-z)r_{zzz}-2r_{zz}-{1\over 8r^3}-{3(y-z)\over 8r^4}r_z+{3(y-z)^2\over 8r^4}r_{zz}-{3(y-z)\over 4r^4}r_z$$$$-{3(y-z)^2\over 2r^5}(r_z)^2+{3(y-z)^2\over 16r^4}r_{zz}-{(y-z)^3\over 16r^4}r_{zzz}+{(y-z)^3\over 4r^5}r_{zz}r_z+{(y-z)^3\over 2r^5}r_{zz}r_z$$$$-{3(y-z)^2\over 4r^5}(r_z)^2-{5(y-z)^3\over 4r^6}(r_z)^3+{2x^2\over (y-z)^3}r_z+{x^2\over (y-z)^2}r_{zz}+{6rx^2\over (y-z)^4}+{2x^2\over (y-z)^3}r_z$$$$+{x^2\over (y-z)^2}r_{zz}+{x^2\over (y-z)}r_{zzz} +{x^2\over (y-z)^2}r_{zz}+{2x^2\over (y-z)^3}r_z.\eqno{(A.13)}$$
Using $r=-1,$ $r_z=0$, $r_{zz}=1/4$ and $y-z=2x^{1/2}$, (A.13) becomes
$$\phi_{zzz}=-{1\over 4}+2x^{1/2}r_{zzz}-{1\over 2}+{1\over 8}+{3x\over 8}+{3x\over 16}-{x^{3/2}\over 2}r_{zzz}
+{x\over 16}-{3\over 8} +{x\over 16}+{x^{3/2}\over 2}r_{zzz}+{x\over 16}$$
$$=2x^{1/2}r_{zzz}-1+{3x\over 4}.\eqno{(A.14)}$$
Substituting $r_{zzz}$ from (A.6), this says $\phi_{zzz}=-1/4$, and that makes $$\Phi^{sp}_{zzz}(t,x,y,-1,y-2\sqrt x)=0.\eqno{(A.15)}$$
That is perhaps surprising, but it will make $w(x,y,t;k)$ of the same order as the incoming beam.

The imaginary part of $\Phi^{sp}_{zzzz}$ is $3/4$. Since this is positive, the  real part of $\Phi^{sp}_{zzzz}$ is essentially irrelevant. However, for completeness we compute it here.  Since we only need $\phi_{zzzz}(x,y,y-2\sqrt x)$, when we differentiate (A.13), we can neglect all terms multiplied by $r_z$. That gives (using no other simplifications)
$$\phi_{zzzz}=-4r_{zzz}+(y-z)r_{zzzz}-{3(y-z)\over 8r^4}r_{zz}-{3(y-z)\over 4r^4}r_{zz}+{3(y-z)^2\over 8r^4}r_{zzz}$$$$-{3(y-z)\over 4r^4}r_{zz}-{3(y-z)\over 8r^4}r_{zz}+{3(y-z)^2\over 16r^4}r_{zzz}+{3(y-z)^2\over 16r^4}r_{zzz}-{(y-z)^3\over 16r^4}r_{zzzz}$$$$+{3(y-z)^3\over 4r^5}(r_{zz})^2+{4x^2\over (y-z)^3}r_{zz}+{x^2\over (y-z)^2}r_{zzz}+{24rx^2\over (y-z)^5}+{2x^2\over (y-z)^3}r_{zz}+{x^2\over (y-z)^2}r_{zzz}\hskip1.5in$$$$+{2x^2\over (y-z)^3}r_{zz}+{x^2\over (y-z)}r_{zzzz}+{2x^2\over (y-z)^2}r_{zzz}+{2x^2\over (y-z)^3}r_{zz}+{2x^2\over (y-z)^3}r_{zz}$$
This can now be further simplified using $r=-1$, $r_{zz}=1/4$ and $2x^{1/2}=y-z$. That gives
$$\phi_{zzzz}=2x^{1/2}r_{zzzz}+(-4+4x)r_{zzz}-{3\over 8}x^{3/2}-{3\over 4}x^{1/2}-{3\over 4}x^{-1/2}.\eqno{(A.16)}$$
Substituting for $r_{zzzz}$ and $r_{zzz}$ from (A.8), we have
$$\phi_{zzzz}={3\over 8}(x^{1/2}-x^{-1/2}), \eqno{(A.17)}$$
and  $\Phi^{sp}_{zzzz}=\nu\phi_{zzzz}+3i/4$.
\vskip.2in
\centerline{\bf Appendix III}
\vskip.2in
Using (3) and changing $y$ to $z$, $\hat f(\eta, \tau)$  is given by
$$\hat f(\eta, \tau)=\int_{\Bbb R^2}a(z)e^{-i(\eta z+\tau t -k\psi(0,z,t))}dtdz,$$
where
$$\psi(0,z,t)=-{z^3\over 8} -(t-t(z))+{1\over 2}(-{z^2\over 4},t-t(z))\cdot M(z)(-{z^2\over 4},t-t(z)).$$
Using the standard notation for the entries in the symmetric matrix $M$ and substituting $t=s+t(z)$, we have
$$\hat f(\eta, \tau)=\int_{\Bbb R^2}a(z)e^{-i\tilde \psi(s,z)}dsdz,$$
where
$$\tilde \psi(s,z)=\tau t(z)+\eta z +{k\over 8}z^3-{k\over 32}z^4M_{11}(z)+(\tau +k-{k\over 4}z^2M_{12}(z))s-{k\over 32}M_{22}s^2.$$
Evaluating the gaussian integral in $s$ and substituting $\tau=k\nu$, $\eta=k\mu$
$$\hat f(k\mu,k\mu)=\int_{\Bbb R}\big({2\pi\over-ikM_{22}(z)}\big)^{1/2}a(z)e^{-ik\rho(z,\mu,\nu)}dz,$$
where
$$\rho=\nu t(z)+\mu z +{1\over 8}z^3-{1\over 32}z^4M_{11}(z)+(2M_{22}(z))^{-1}(\nu+1+{1\over 4}z^2M_{12}(z))^2\eqno{(A.18)}$$

All of the functions of $z$ in (A.18) are carried through the computations without change until the steepest descent computation in \S 5. We need to compare (A.18) with (4). In (4) $-{1\over 32}z^4M_{11}(z)$ is replaced by $-{i\over 32}z^4$. In view of (3) ${1\over 32}z^4M_{11}(z)={i\over 32}z^4+O(z^5)$. In (4) $(2M_{22}(z))^{-1}(\nu+1+{1\over 4}z^2M_{12}(z))^2$ is replaced by $-{i\over 2}(\nu+1)^2$. To see the effect of that consider (18): the function $C$ is unchanged, and the stationary point in $\nu$ is now
$$\nu=-1-{1\over 4}z^2M_{12}(z)-2M_{22}(z)C.$$ Since $M_{12}=O(z)$ and $C=O(z^3)$, we have $\nu=-1+O(z^3)$ as before. All of this -- including the correction from ${1\over 32}z^4M_{11}(z)$ -- changes the final exponent $i(B-C)-C^2/2$ in (18), to $i(B-C)-C^2/2+O(z^5)$. The $O(z^5)$ becomes negligible in \S 6. Similarly, $\big({2\pi\over-iM_{22}(z)}\big)^{1/2}a(z)=(2\pi)^{1/2}+O(z)$, so that the change in the amplitude in (18) is negligible as well.
\newpage

\noindent {\bf References}
\vskip.2in
\noindent [F] Friedlander, F.G., The wave front set of the solution of a simple initial-boundary value problem with glancing rays, Math. Proc. Cambridge Phil. Soc. {\bf 79}(1976), 145-159.
\vskip.1in
\noindent [H] H\"ormander, L,. The
Analysis of Linear Partial Differential Operators, I-IV,
Springer-Verlag, Berlin (1985)
\vskip.1in
\noindent [L] Lebeau, G. Propagation des ondes dans les variet\'es \`a coins, Ann. Scient. \'Ecole Norm. Sup. 30 (1997), 429-497.
\vskip.1in
\noindent [MT] Melrose, R. and Taylor, M.E., Near peak scattering and the corrected Kirchhoff approximation for convex obstacles, Advances in Math. 55 (1985), 242-315.
\vskip.1in
\noindent [R] Ralston, J.,  Gliding beams and whispering gallery modes,

\noindent http://www.math.ucla.edu/~ralston/pub/Glide2.pdf
\vskip.1in
\noindent [R1] Ralston, J., Construction of Gaussian Beam Phases for Reduced Eichonals,
 http://www.math.ucla.edu/~ralston/pub/reducedbeam.pdf

\vskip.1in
\noindent [Sj] Sj\"ostrand, J. Propagation of analytic singularities for second order Dirichlet problems, CPDE {\bf 5} (1980), 41-93.
\vskip.1in
\noindent [T] Taylor, M.E., Grazing rays and reflection of singularities of solutions to wave equations, CPAM {\bf 29} (1976), 1-38.
\vskip.1in
\noindent [Ti] Tiruviluamala, N. On the Passage of Gaussian Beams through Cusps in Ray Paths, http://www.math.ucla.edu/~ralston/pub/thesis.pdf

\end